\documentclass[letterpaper,11pt]{article}
\usepackage{amsmath, amsthm, amssymb}    
\usepackage{bridges}			
\usepackage{graphicx}			
\usepackage[colorlinks=true, urlcolor=blue, citecolor=black, linkcolor=black]{hyperref}  
\usepackage{subcaption}			

\usepackage{algpseudocode}

\newcommand{\RR}{\mathbb{R}}
\newcommand{\ZZ}{\mathbb{Z}}
\newcommand{\FF}{\mathbb{F}}
\newcommand{\CC}{\mathbb{C}}

\renewcommand{\SS}{\mathbb{S}}

\title{Elliptic Curves  and the Hopf Fibration}

\author{Nadir Hajouji\textsuperscript{1} and Steve Trettel\textsuperscript{2}
\vspace{10pt}\\
\textsuperscript{1}Boston, MA, USA; nhajouji@gmail.com\\
\textsuperscript{2}University of San Francisco, CA, USA; strettel@usfca.edu} %
\date{}					

\begin{document}

\maketitle

\thispagestyle{empty}

\begin{abstract}
By combining tools from different areas of mathematics,
we obtain 3D visualizations of elliptic curves over different fields that faithfully capture the underlying algebra and geometry.

\end{abstract}

\vspace{-3mm}
 \section*{Introduction}

This paper is about visualizing elliptic curves.
If you've never heard of an elliptic curve,
here's all you need to know: 1) They are not ellipses, 2) you might not recognize them as curves,
3) they are \emph{incredibly} interesting, appearing across mathematics, cryptography, and physics and 4) they are \emph{notoriously} hard to explain to non mathematicians. Although they can be defined in simple terms, the modern definition does very little to explain their ubiquity, or exhibit their beauty.

We want to show you what elliptic curves really look like. Using tools from complex analysis, algebra, number theory and topology, we bring these curves into view as 3D renders, and produce a gallery showcasing their diversity.
We aim to make elliptic curves accessible to newcomers, and to reveal their aesthetic beauty to a mathematical audience, who met them first through symbols on a page.

\begin{figure}[h!tbp]
	\centering
	\includegraphics[width=0.9\textwidth]{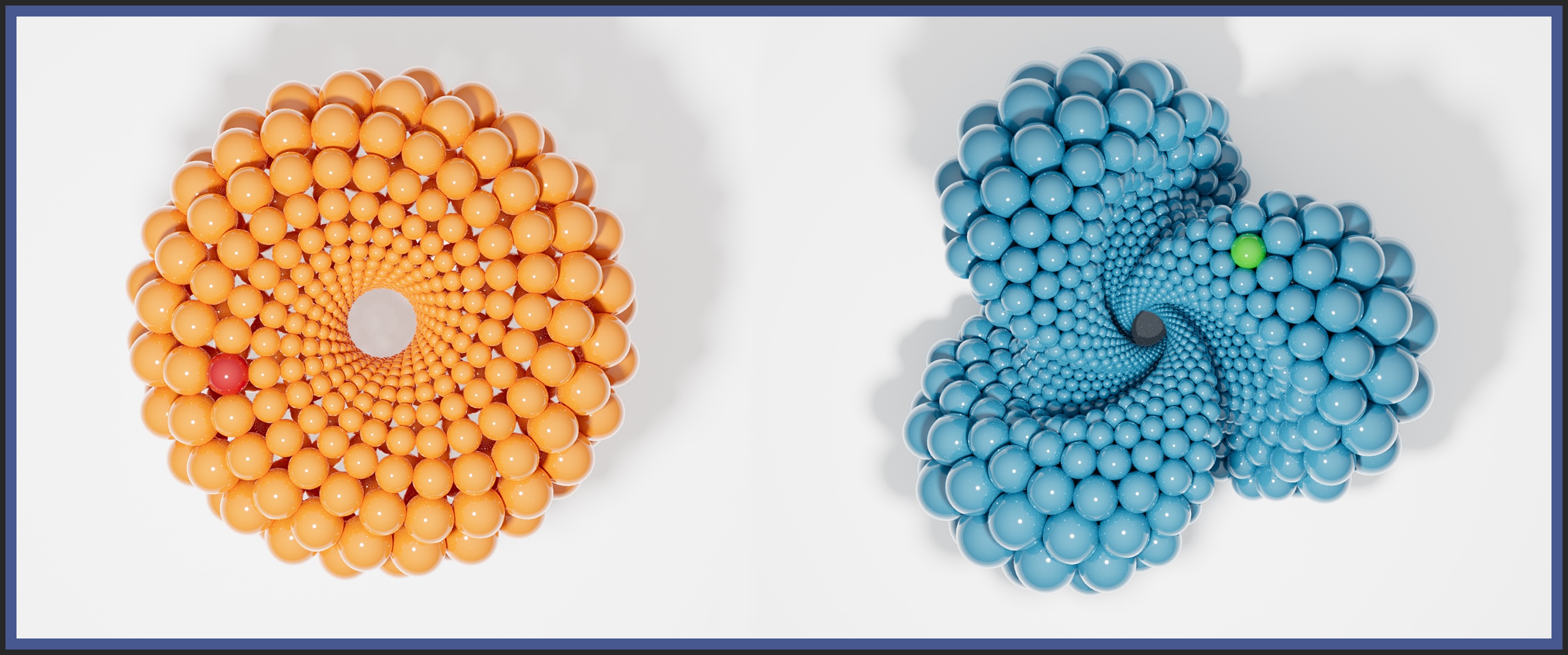}
    \caption{The elliptic curves $y^2=x^3+3x$ over $\FF_{625}$ (orange), and $y^2=x^3+3$ over $\FF_{2401}$ (blue).}
	\label{fig:bestPic}
\end{figure}

\section*{What is an Elliptic Curve?}

\emph{Everyone knows what a curve is, until they have studied enough mathematics to become
confused.\\
— F. Klein\footnote{The idea to use quotations to start our sections, as well as the quotes themselves, came from \cite{RisingSea}.}}

\vspace{2mm}

When a classical geometer says the word ``curve," they mean \emph{something you can draw on a (possibly infinite) piece of paper with a pen}. The unifying concept is ``one dimensionality": each point on the curve has a left and a right, each moment during its drawing a before and an after. 
In modern mathematics, the word \emph{curve} may refer to any one dimensional object \textemdash anything that can described by a single variable.
The invention of abstract algebra has allowed us the freedom to consider variables over many wild and wonderful number systems, and opened our eyes to a more varied collection of curves.
Real curves look familiar, like segments of the number line,
but there are also \emph{complex curves}, which trace out something we might more readily recognize as a surface,
and curves defined over \emph{``finite number systems"},
which look like a cloud of isolated points.
All of these are equally curves in the eyes of an algebraic geometer: a ``single variable object" across different mathematical universes.  It is in this modern, tolerant sense that elliptic curves are curves; seeking a unifying perspective for this varied family is one of our mathematical and artistic motivations.

\begin{figure}[h!tbp]
	\centering
    \includegraphics[width=0.9\textwidth]{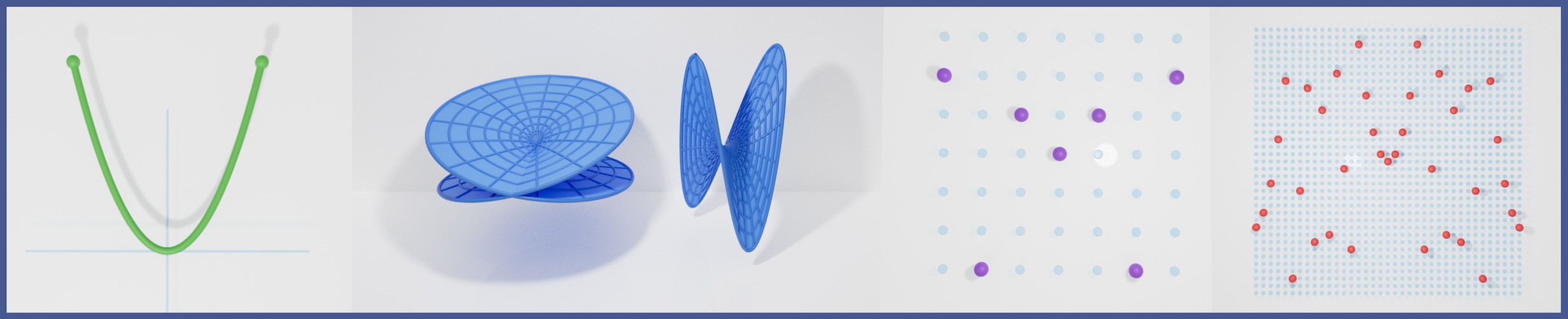}
	\caption{The curve $y=x^2$ over $\RR$ (green), over $\CC$ (blue: two views, projecting to three dimensions by deleting the imaginary component of $x$ or $y$), and over finite fields: $\FF_7$ (purple), $\FF_{37}$ (red).
}
	\label{fig:abstract-curve}
\end{figure}

\subsection*{Perspectives on the Circle}
Circles are perhaps the most familiar examples of mathematical curves,
and we will see that elliptic curves are a natural generalization of circles.
To make the similarity between the two more apparent,
we start by answering the following question \textemdash \emph{how do you describe a circle in algebraic language?}
We will give 3 distinct ways of describing the unit circle, and subsequently show that elliptic curves can be defined by tweaking any of those 3 descriptions.
The first description, due to Pythagoras and Descartes,
is via an equation: the unit circle is the set of solutions $x,y\in\RR$ to  $x^2+y^2=1$.
The second description is via a parametrization: the solutions to this equation coincide with the values of $(\cos \theta, \sin \theta)$ as $\theta$ varies across the set of real numbers.
Of course, we don't need to use every real number \textemdash replacing $\theta$ by $\theta+2\pi m$,
where $m$ is any integer,
does not change the point, so we really only need $\theta$ between 0 and $2\pi$, or between $-\pi$ and $\pi$.
This leads to our third description:
points $\theta$ on the circle are like real numbers, except subject to the condition that they ``do not change" if we add integer multiples of $2\pi$.
That is, the circle is made out of ``ill defined real numbers": or numbers whose values only makes sense up to this $\pm2\pi n$ ambiguity.  Mathematicians write such numbers as $\RR/2\pi \ZZ$, and more generally $G/H$ (read $G$ \emph{mod} $H$) for elements of $G$ which are ambiguous up to elements of $H$.

\begin{figure}[h!tbp]
	\centering
    \includegraphics[width=0.9\textwidth]{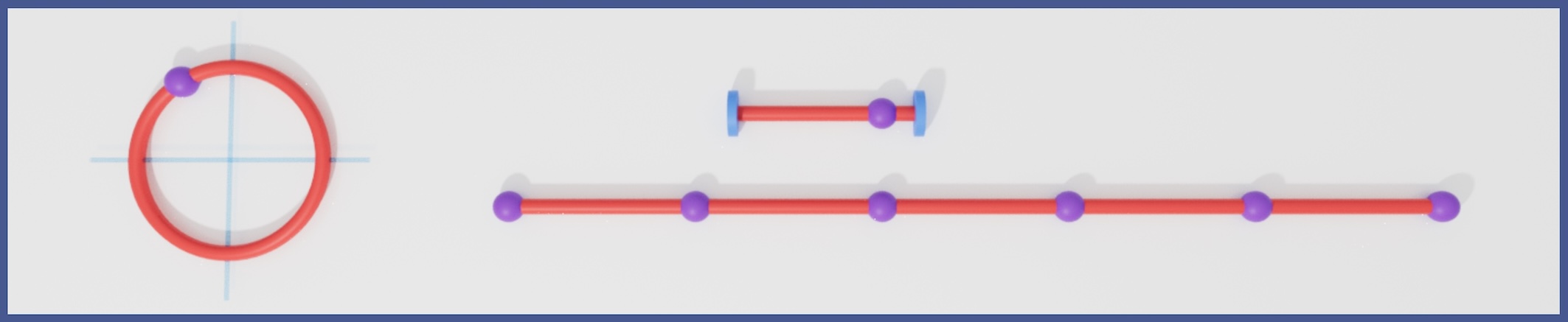}
	\caption{A circle three ways: (1,2) an implicit / parametric plane curve in $\RR^2$ (left); (3) an ill-defined real number in $\RR/2\pi\ZZ$ as a periodic point in $\RR$ or element of a fundamental domain (right).}
	\label{fig:circle}
\end{figure}

\subsection*{Elliptic Curves}

Elliptic curves are one of the simplest shapes, other than circles,
that can be described in each of these three ways.
Our starting point will be a generalization of (3):
we replace the real numbers by complex numbers, the set of integer multiples of $2\pi$
by a lattice $\Lambda = \omega_1 \ZZ \oplus \omega_2 \ZZ$,
and define an elliptic curve to be the set of ''ill-defined complex numbers"
$\CC/\Lambda$.
Topologically such numbers form a \emph{torus}: a quite natural generalization of the circle indeed! 
From this starting point, we can build back the other perspectives of our triad.  The analog of $\cos \theta, \sin \theta$ for an elliptic curve are the functions $\wp_\Lambda, \wp_\Lambda'$, known as the \emph{Weierstrass $\wp$-function} $ \wp_\Lambda(z) = \frac{1}{z^2} + \sum_{\lambda - \Lambda - \{0\}} \frac{1}{(z-\lambda)^2} - \frac{1}{\lambda^2}$ and  its derivative $\wp_\Lambda'(z) = -2\sum_{\lambda \in \Lambda } \frac{1}{(z-\lambda)^3}$ 
They satisfy an equation of the form:
\begin{equation}
    \left( \frac{\wp_\Lambda'(z)}{2} \right)^2 = \wp_\Lambda(z)^3+\left(15 \sum_{\lambda \in\Lambda - \{0\}} \frac{1}{\lambda^{4}} \right) \wp_\Lambda(z)+\left(35 \sum_{\lambda \in\Lambda - \{0\}} \frac{1}{\lambda^{6}}\right)
\end{equation}
Calling the infinite sums appearing in this differential equation $f_\Lambda$ and $g_\Lambda$ respectively, we see that $(x,y)=(\wp, \frac{1}{2}\wp^\prime)$ satisfy the algebraic equation  $y^2 = x^3 + f_\Lambda x + g_\Lambda$.
This is called the equation of the elliptic curve.

The appearance of these constants $f_\Lambda,g_\Lambda$ reveals an important fact: unlike circles, elliptic curves come in \emph{different flavors}. While there is essentially only one circle (any two are equivalent up to scaling), there are many distinct lattices in the complex plane that are not related by scaling or rotation.
The corresponding elliptic curves inherit distinct personalities from their lattice progenitors. This additional source of variety is part of what makes elliptic curves such a compelling subject for visual art.

\section*{Visualizing Classical Elliptic Curves}

\emph{We all know that Art is not truth. Art is a lie that makes us realize truth, at least the truth that is given us to understand. The artist must know the manner whereby to convince
others of the truthfulness of his lies.\\
— P. Picasso }

\begin{figure}[h!tbp]
	\centering
    \includegraphics[width=0.9\textwidth]{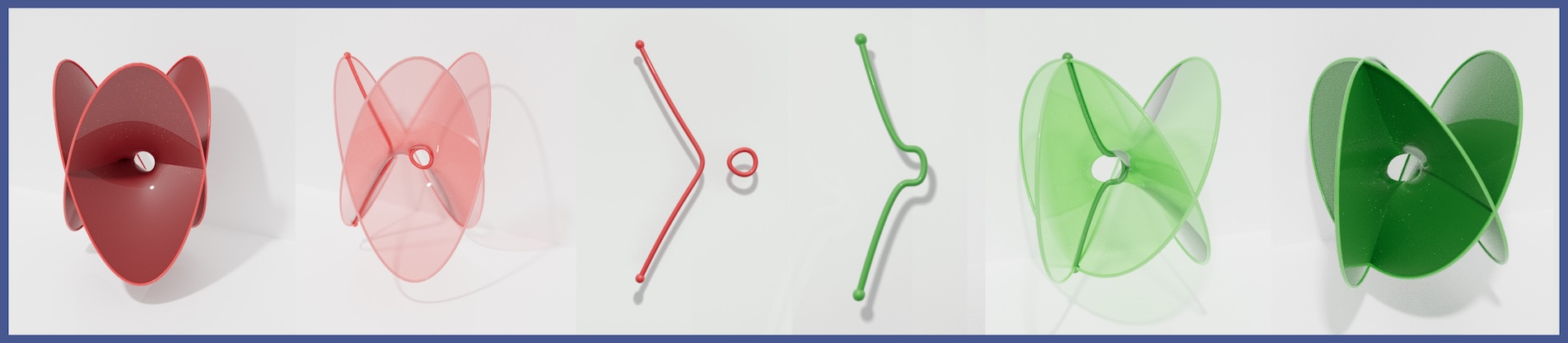}
	\caption{Equation-based illustrations of elliptic curves, for the curves $y^2 = x^3 -x$ (red) and $y^2 = x^3 + 1$ (green). These illustrations obscure the underlying torus geometry.}
	\label{fig:weierstrass}
\end{figure}

\vspace{-3mm}

Having the equations defining elliptic curves readily available, it is tempting to try plotting the solutions, much as one does to visualize the circle.  
There's just one small problem: $x$ and $y$ are complex, so the points $(x,y)$ live in a \emph{four dimensional space}.  
For our three dimensional brains, one must then either throw away a dimension
(introducing artificial self-intersections), or discard most of the shape, plotting a two dimensional \emph{real slice}.  Such images can be seen in Figure \ref{fig:weierstrass}.

\begin{figure}[h!tbp]
	\centering
    \includegraphics[width=0.9\textwidth]{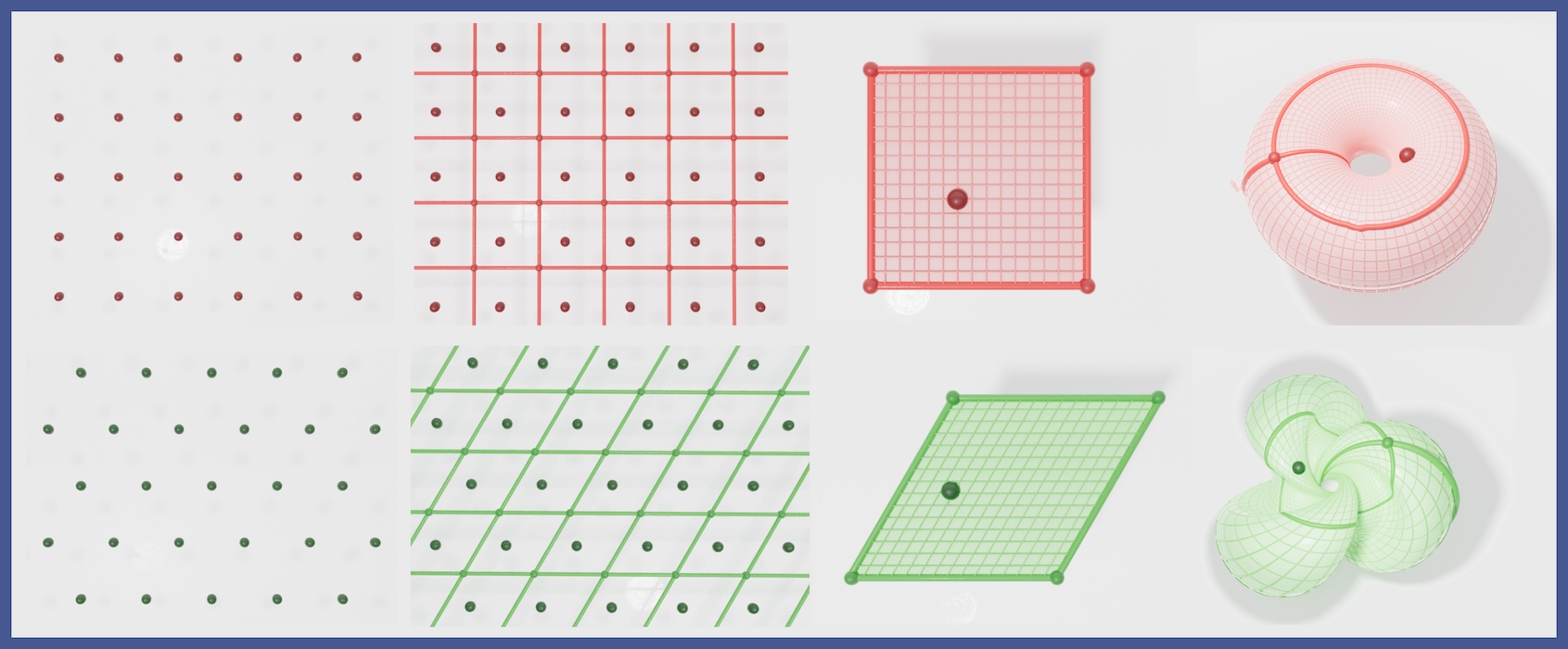}
	\caption{A point on an elliptic curve: three views (1) as lattice (2) in fundamental domain (3) rolled up. The same elliptic curves as the previous illustration.}
	\label{fig:rollup}
\end{figure}

A goal of this paper is to present an alternative approach to visualizing elliptic curves, using their description as quotients $\CC/\Lambda$. This is summarized in Figure \ref{fig:rollup}: beginning with an infinitely repeating grid, one may produce a finite visualization by 
selecting a fundamental domain, although this comes at a cost:
the edges of the chosen domain introduce arbitrary discontinuities, sacrificing symmetry and visual elegance. A more compelling alternative to these flattened images would seamlessly ``roll up" this repetition onto a smooth toroidal surface, analogous to wrapping a repeating real line around a circle.  Constrained by our biology and the local laws of physics, we seek a means of doing this in three dimensions (as opposed to the mathematically natural four-dimensional candidate $(\wp,\wp^\prime)$).

\vspace{-5mm}

\subsection*{Hopf to the Rescue}

While constructing such a ''roll up map" into three dimensions, we must be mindful not to destroy the distinct personalities imbued by their different underlying lattices. 
Complex geometry tells us any allowable map $\CC/\Lambda\to\RR^3$ must be angle-preserving (conformal): a strong constraint which means we can't just roll $\CC/\Lambda$ up as a torus of revolution, for instance.
One natural source of conformal maps is \emph{isometries}, so one might hope to seek an isometric embedding $\CC/\Lambda\to\RR^3$.  Unfortunately a classic argument in differential geometry says we cannot do this: there are no flat tori in $\RR^3$ and thus no isometries from $\CC/\Lambda$ into $\RR^3$.
However one can save the idea with a small tweak: since stereographic projection is conformal, it would be enough to find flat tori into the 3-sphere, and then map them into $\RR^3$ by a composition of conformal maps!

\begin{figure}[h!tbp]
	\centering
    \includegraphics[width=0.9\textwidth]{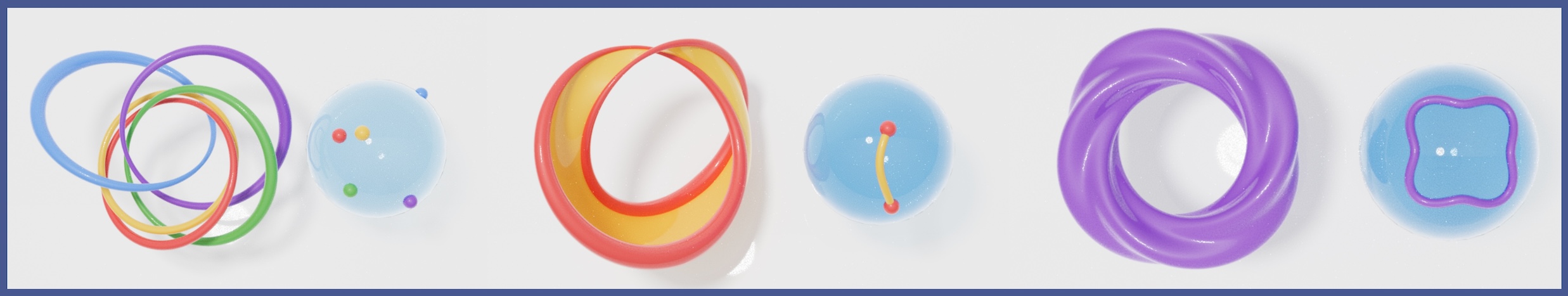}
	\caption{Point preimages under the Hopf fibration are circles, so the preimages of curves are annuli or tori.}
	\label{}
\end{figure}

Amazingly, such flat tori exist and were were described by Ulrich Pinkall in \cite{pinkall}. 
The key tool is the Hopf fibration, a topologically nontrivial map $\eta\colon \SS^3 \to \SS^2$ defined by $\eta(z, w) = z/w$ for unit vectors $(z, w) \in \CC^2$, with $z/w$ interpreted as a point on the Riemann sphere $\CC \cup {\infty}$. This simple formula conceals beautiful geometric structure: each point’s preimage is a circle, so an arc lifts to an annulus, and a closed curve to a torus. 
A straightforward calculation shows that these tori are flat, setting up a correspondence between tori $\CC/\Lambda$, and closed curves on the 2 sphere.
 Pinkall makes this explicit, showing the preimage $\eta^{-1}(C)$ of a simple closed curve $C$ on $\mathbb{S}^2$ of length $L$ enclosing area $A < 2\pi$ is isometric to $\mathbb{C} / \Lambda$, for $\Lambda = 2\pi \ZZ\oplus (\tfrac{A}{2}+i\tfrac{L}{2}) \ZZ$.  
Thus, to realize an elliptic curve  $\CC/\Lambda$,
we need only find a loop on the sphere with the appropriate $A$ and $L$.

\begin{figure}[h!tbp]
	\centering
	\includegraphics[width=0.9\textwidth]{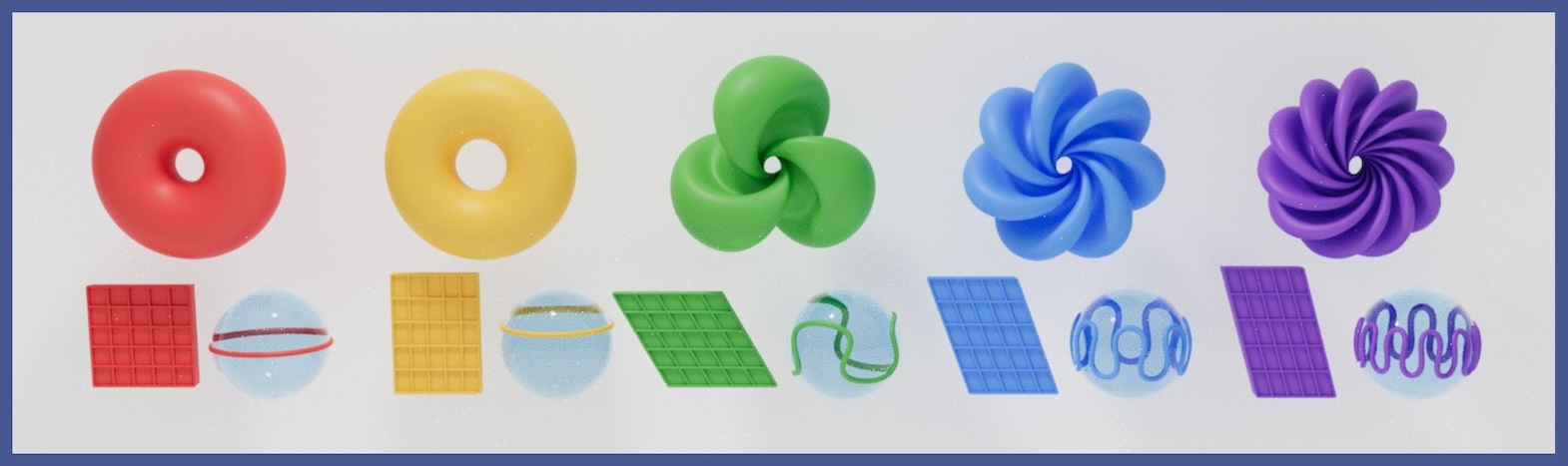}
	\caption{Lattices generated by $1$ and $\tau\in\{i,i\sqrt{2},\tfrac{1+i\sqrt{3}}{2},\tfrac{1+i\sqrt{7}}{2},\tfrac{1+i\sqrt{11}}{2}\}$, curves on the sphere with the corresponding lengths and areas, and the resulting embedded elliptic curves in $\RR^3$.}
	\label{fig:HopfTori}
\end{figure}

Fixing area and length still leaves a generically infinite dimensional space of possible curves, which provide  substantial artistic freedom in constructing embeddings. 
After choosing an appropriate curve $C$, the real computation begins: we need the isometry $\CC/\Lambda\to\eta^{-1}(C)$.
We sketch this below, slightly generalizing the techniques of \cite{banchoff}.
Parameterizing $C$ by $(\theta(t),\phi(t))$ and the Hopf fiber above $(\theta,\phi)\in\SS^2$ as $H_{(\theta,\phi)}(s)=(e^{i(\theta+s)}\sin\tfrac{\phi}{2},e^{is}\cos\tfrac{\phi}{2})$ we can parameterize the torus $\eta^{-1}(C)$ by $s+it\mapsto H_{c(t)}(s)$ in $\SS^3$.  While this initial parameterization is not an isometry, we modify it into one by accounting the curvature of the Hopf fibration (as a circle bundle over $\SS^2$), and then compose with stereographic projection into $\RR^3$.  To save others the trouble of computing this, we give the resulting map $\CC/\Lambda\to\RR^3$ below:

{\begin{algorithmic}[1]
\State {\bfseries Given} $s+it\in\CC$ with $t<L/2$ and a curve $(\theta(x),\phi(x))$ on $\SS^2$:
\State Numerically find $v$ such that $L(v)=\int_0^v (\theta^\prime(x)^2\sin^2\phi(x)+\phi^\prime(x)^2)^{1/2}\,dx=2t$
\State Compute $\theta=\theta(v)$, $\phi=\phi(v)$ and $f=\int_0^v \sin\left(\phi(x)/2\right)\theta^\prime(x)\,dx$
\State Compute $h=H_{(\theta,\phi)}(s-f)=(e^{i (\theta+s-f)}\sin\tfrac{\phi}{2} ,e^{i(s-f)}\cos\tfrac{\phi}{2})$
\State {\bfseries Return} stereographic projection $\sigma(h)$ for $\sigma\colon (x,y,z,w)\mapsto (x,y,z)/(1-w)$
\end{algorithmic}}

\section*{Elliptic Curves over Other Number Systems}

\emph{Algebra is the offer made by the devil to the mathematician. The devil says: I will give
you this powerful machine, it will answer any question you like. All you need to do is give
me your soul: give up geometry and you will have this marvelous machine.
— M. Atiyah}

Taking the modern mathematician's notion of curve seriously, we have become well acquainted with many donut-shaped elliptic curves.  But these are just one lineage in the vast elliptic curve family: to meet others, we must learn how to swap $\CC$ for other number systems. This requires using a different part of our triad: we will need the equation of the elliptic curve. 
Given a number system $k$, an \emph{elliptic curve over $k$} is the solutions $x,y\in k$ to an equation
of the form $y^2=x^3+fx+g$.
We might then hope to visualize these new curves by plotting their solutions in $(x,y)\in k^2$ (at least when we can reasonably visualize $k$).  We've actually already seen such plots for $\RR$ in Figure \ref{fig:weierstrass}, and plots over finite number systems look like scatterings of dots akin to Figure \ref{fig:abstract-curve}. 
This may give the impression that the family resemblance between the various lineages of elliptic curves ends with their equation.  But that is not the case: despite vast algebraic differences, a deep geometric unity lurks beneath the surface.  Looking over many fields, we'll find that classical elliptic curves survive in the background like ghosts, subtly influencing their structure.

Our goal in this section is to realize this unified view, and render images across the elliptic curve family that highlight familiar toroidal shapes and lattice-like symmetries.
To do so we employ Galois theory \textemdash an advanced toolkit from algebra that was developed by one of the romantic heroes of mathematics.\footnote{Indeed, Galois came up with the tools we need when he was a mere 16 yeas old. His genius was lost to the world after his untimely death in a duel at the tender age of 20.}
Galois teaches us: \emph{the key to understanding algebra over small number systems lies in understanding symmetries of larger number systems, like the complex numbers}.

\begin{figure}[h!tbp]
	\centering
    \includegraphics[width=0.9\textwidth]{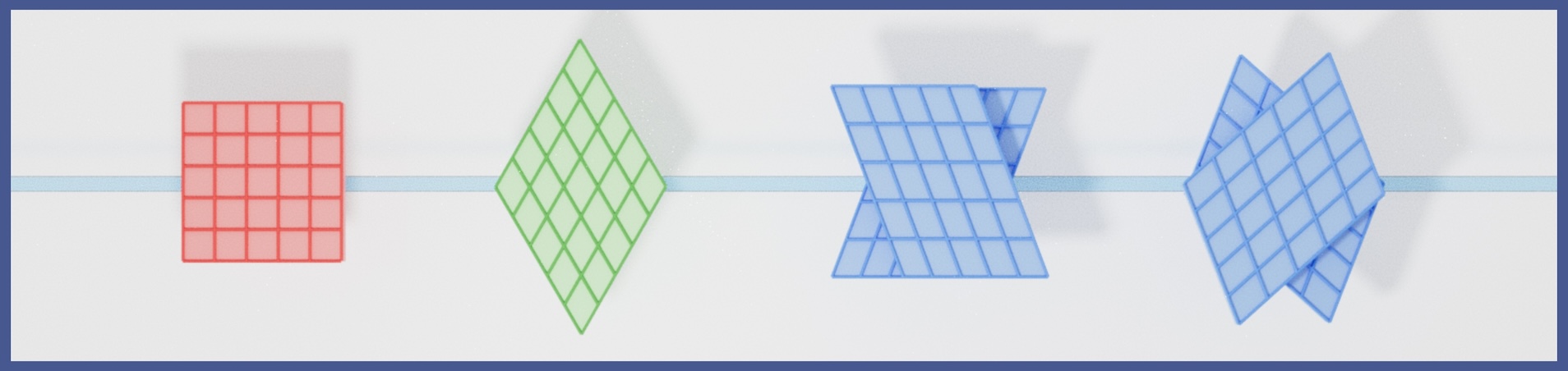}
	\caption{A pair of ``real" lattices (red, green). ``Non-real" (i.e. complex) lattices (blue) have no reflection symmetry, in any orientation.
}
	\label{fig:realLattices}
\end{figure}

\begin{figure}[h!tbp]
	\centering
	\includegraphics[width=0.9\textwidth]{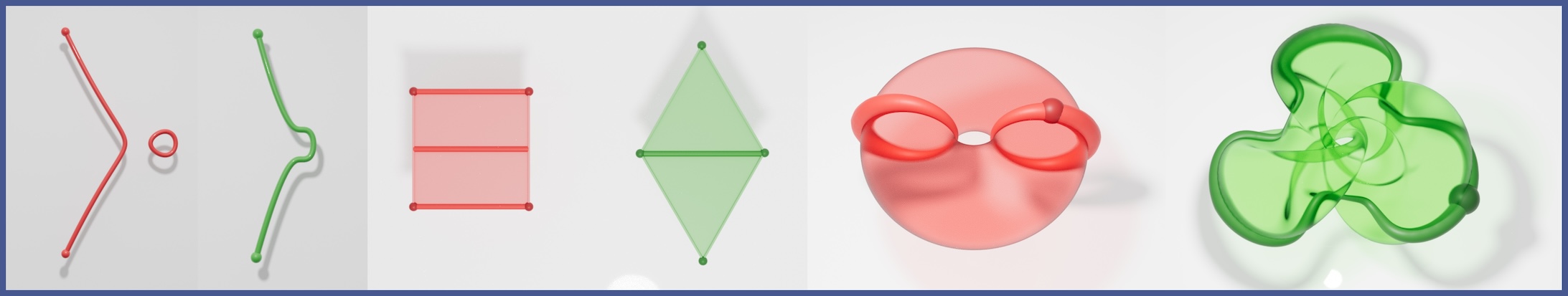}
	\caption{The real elliptic curves $y^2=x^3+3x$ (red) and $y^2=x^3+1$ (green), as standardly depicted as solutions in $\RR^2$ (left), and as subsets of a classical elliptic curve $\CC/\Lambda$ (center, right). The standard illustrations close up via a point at infinity, highlighted as a finite point when embedded in $\CC/\Lambda$. }
	\label{fig:realCurves}
\end{figure}

We explain this cryptic decree via example, using Galois' insight to find real elliptic curves among the complex tori we already understand.
Our first step is to realize one can talk about $\RR$ inside of $\CC$ using the language of symmetry: real numbers are precisely those that are \emph{unchanged} by complex conjugation.
With some work \cite{ECwebsite}\cite{SilvermanAdv}, we can show that this gives a way to discover which complex curves have a real elliptic curve hiding inside them \textemdash the ones whose lattices are unchanged by reflection in the $x$-axis (the geometric realization of this complex conjugation symmetry), see Figure \ref{fig:realLattices}.  
And once we've sorted these out, we can further uncover the \emph{points} of the hidden real curve itself by seeing what is fixed by complex conjugation in $\CC/\Lambda$. These include the $x$ axis, as well as points $z$ related to $\bar{z}$ by a lattice element.  
With the points in hand, we can depict the real curve inside its complex counterpart using the \emph{roll-up} map derived from the Hopf fibration, or recover the standard illustration in $\RR^2$ via the parameterization $(\wp_\Lambda, \wp_\Lambda')$, see Figure \ref{fig:realCurves}.

\begin{figure}[h!tbp]
	\centering
	\includegraphics[width=0.9\textwidth]{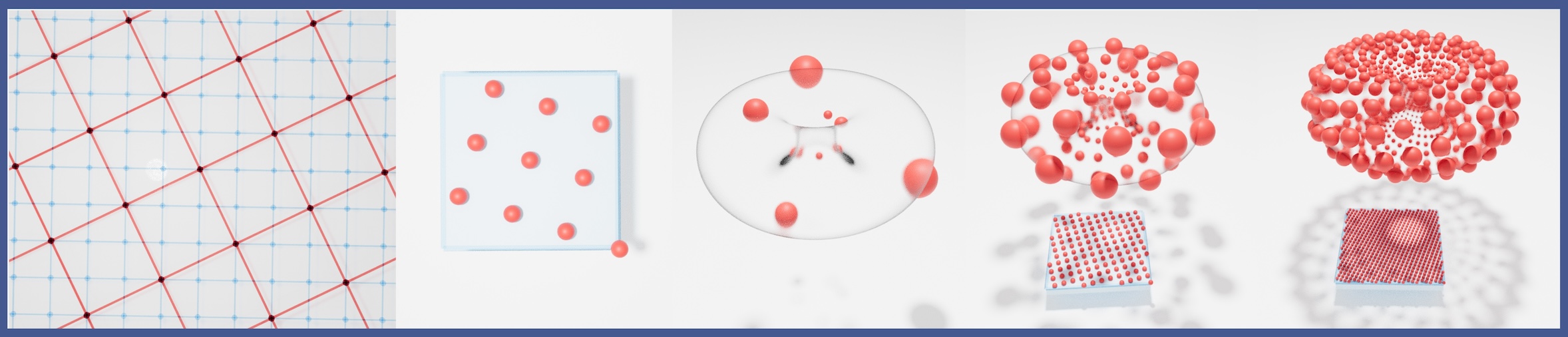}
	\caption{For the elliptic curve $y^2 = x^3+3x\pmod 5$, the Frobenius looks like ''multiplication-by-$(-2+i)$" on the square lattice (far left). We find the $\FF_5$ points by identifying the fixed points of Frobenius (center, on fundamental domain and rolled up), and we similarly we obtain the $\FF_{125}$, $\FF_{625}$ points by identifying fixed points of suitable powers of Frobenius (right).}
	\label{fig:frobenius}
\end{figure}

Galois' philosophy can also be used to visualize elliptic curves over much smaller number systems called \emph{finite fields}.
If these are new, don't worry about the details
\textemdash all you need to know is that the symbol $\FF_q$ denotes a ``number system" that contains precisely $q$ numbers,
including some familiar ones like $0,1,2..$.
To seek out the $\FF_q$ elliptic curves hiding among the complex ones, we will again look for those with a special symmetry: the analog of complex conjugation for $\FF_q$, known as \emph{the Frobenius}.
As a result of some beautiful mathematics \cite{coxprimes}\cite{ECwebsite}\cite{SilvermanAdv}, 
we can show that this map is simple to compute when we represent our curves as $\CC/\Lambda$, acting as ''multiplication by a complex number".
Upon finding a complex curve with this symmetry, we uncover the hidden $\FF_q$ points by seeing what the Frobenius fixes, as in Figure \ref{fig:frobenius}.
This allows us to draw pictures of the $\FF_q$ curve on our familiar complex tori \textemdash even though $\FF_q$ is not a subset of $\CC$!
Furthermore, each lattice-complex number pairing is unique to the elliptic curve,
so we obtain a unique picture for each curve.

\section*{A Gallery of Elliptic Curves over Finite Fields }

\emph{''Should you just be an algebraist or a geometer?” is like saying ''Would you rather be
deaf or blind?”\\
— M. Atiyah.}

At this point in the story,
our training as mathematicians compels us to explain precisely how we're using Galois theory to make our pictures, and how we know our method always works.
We will answer these questions in a later work, and for now, content ourselves with sharing the end result of our labor: new illustrations of elliptic curves over finite fields, living inside the ghost of the complex elliptic curve that once gave birth to it.

\begin{figure}[h!tbp]
	\centering
	\includegraphics[width=0.8\textwidth]{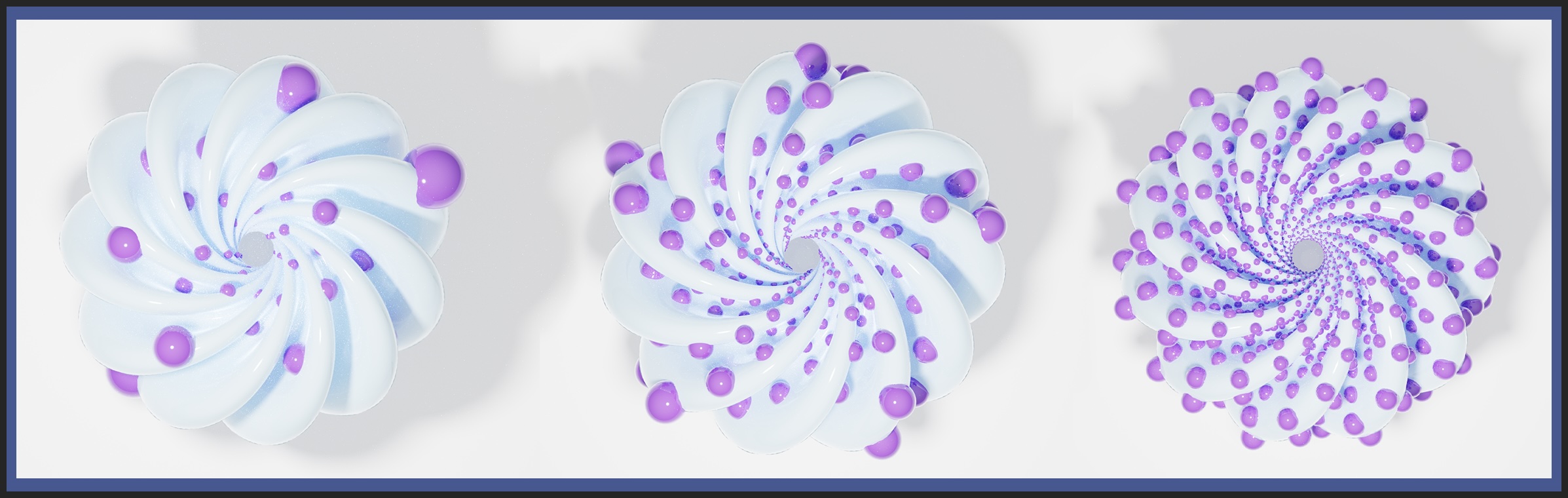}
	\caption{The elliptic curve $y^2=x^3+x+1 \pmod 5$  over $\FF_{125}$ (left), $\FF_{625}$ (middle) and $\FF_{3125}$ (right).}
	\label{fig:D11}
\end{figure}

\begin{figure}[h!tbp]
	\centering
	\includegraphics[width=0.8\textwidth]{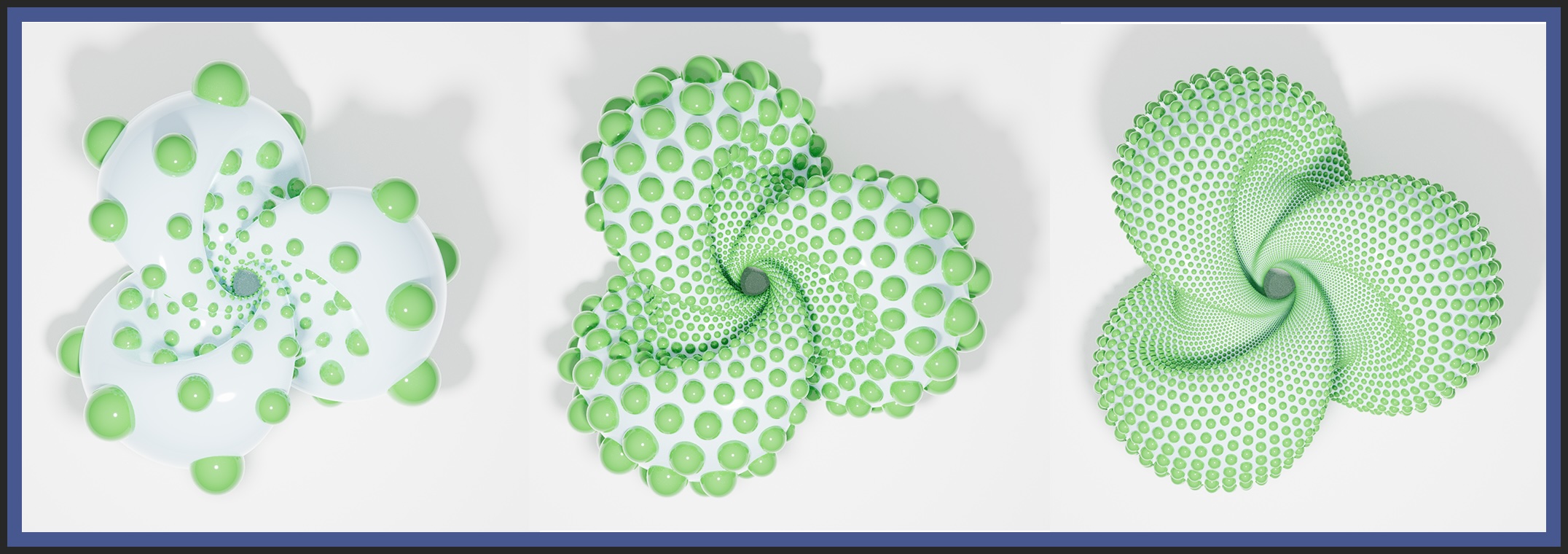}
	\caption{Points of $E\colon \,y^2=x^3+3\pmod 7$ over the fields $\FF_{343}$ (left), $\FF_{2401}$ (center) and $\FF_{16,807}$ (right).  }
	\label{fig:hex-2}
\end{figure}

\begin{figure}[h!tbp]
	\centering
	\includegraphics[width=0.8\textwidth]{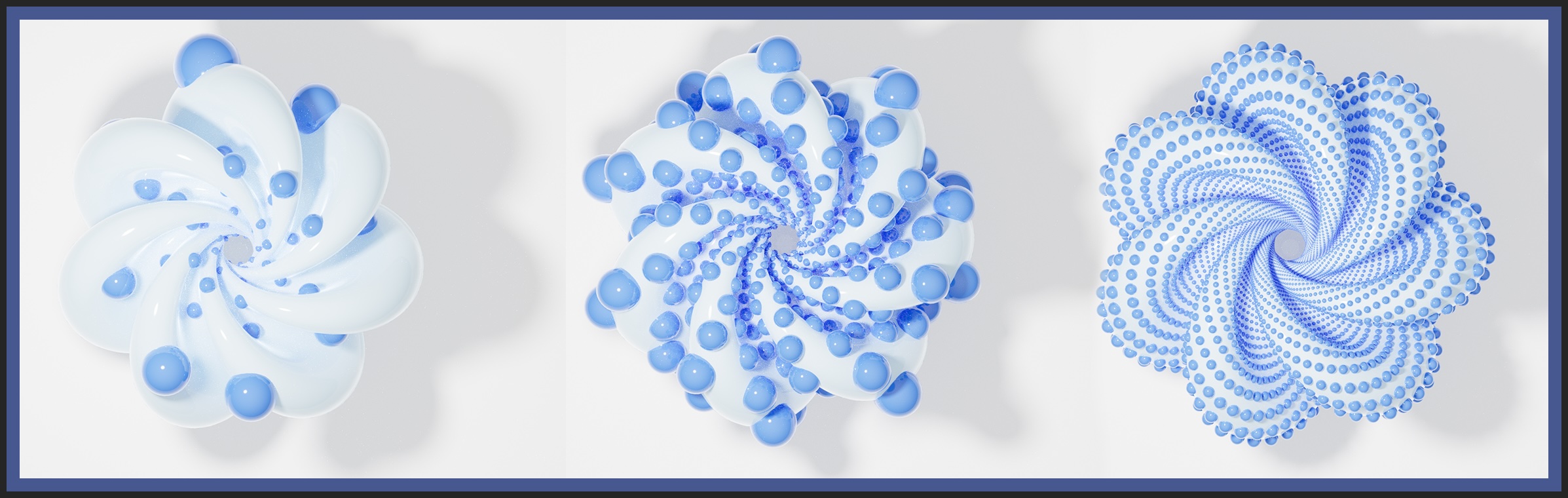}
	\caption{The elliptic curve $y^2=x^3+5x+7\pmod {11}$ over $\FF_{121}$ (left), $\FF_{1331}$ (middle) and $\FF_{14,641}$ (right).}
	\label{fig:D7}
\end{figure}

\begin{figure}[h!tbp]
	\centering
	\includegraphics[width=0.8\textwidth]{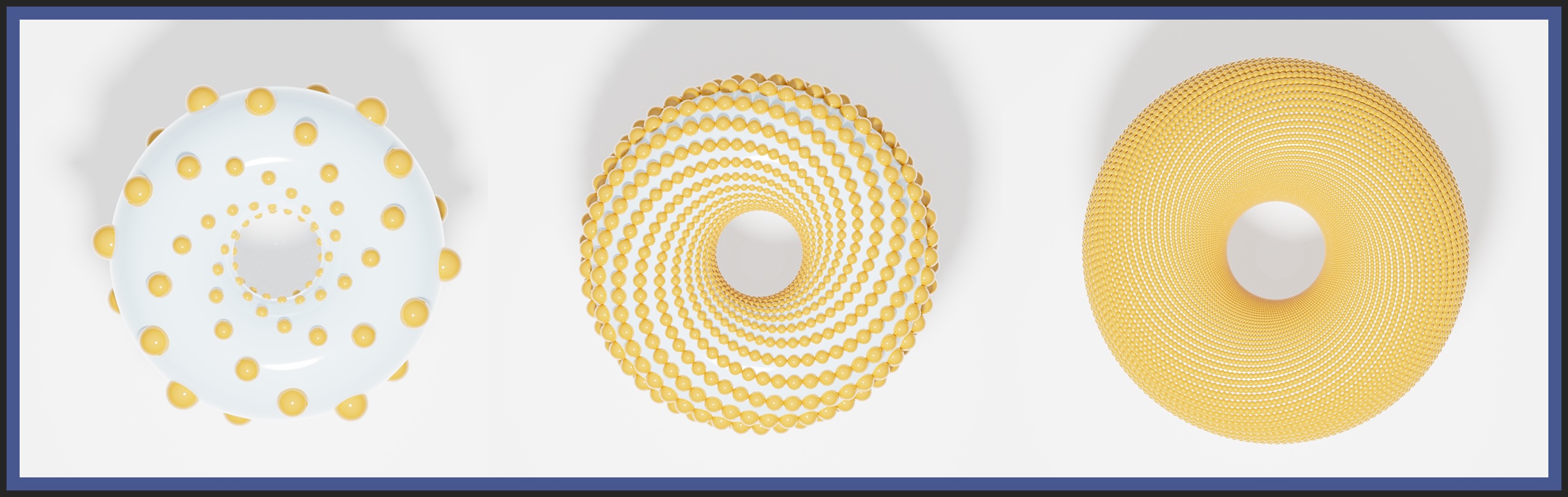}
	\caption{The elliptic curve $y^2=x^3+x+3\pmod {11}$ over $\FF_{121}$ (left), $\FF_{1331}$ (middle) and $\FF_{14,641}$ (right).}
	\label{fig:D8}
\end{figure}
\section*{Conclusion}

In the course of writing this paper, we realized that we would have to answer the following question \textemdash \emph{How do we choose the elliptic curves to include in the paper?} There are infinitely many to choose from, each with its own unique personality,
more than we can ever hope to see in our finite lifetimes,
and we have only 8 pages.
We shared some of our favorites \textemdash but we have a lot more to show!

The true beauty that enthralls mathematicians lies in the patterns that emerge when one looks at \emph{lots} of elliptic curves.
You can think of elliptic curves as musical notes\textemdash we've hopefully convinced you that they ``sound different" from each other, but just like musical notes,
they can only do so much when played in isolation.
We invite you to visit our website \cite{ECwebsite},
where you can see a larger (though always incomplete) gallery and learn about the beautiful mathematics behind our pictures.

\vspace{-0.25cm}
    
{\setlength{\baselineskip}{13pt} 
\raggedright				
\bibliographystyle{bridges}
\bibliography{refs}
} 

\end{document}